\documentclass[12pt]{amsart}
\usepackage{amssymb,amsmath,amsthm,mathtools,amsfonts,times,mathrsfs,multirow}
\usepackage{pgfplots}
\usepackage[hidelinks]{hyperref}

\pgfplotsset{compat=1.15}
\usepackage{mathrsfs,mathdots}
\usetikzlibrary{arrows}
\usepackage{enumerate}
\usepackage{graphicx}
\usepackage{array}
\usepackage{ytableau}
\usepackage{pstricks,pst-node,graphicx}
\usepackage{tikz,varwidth}
\usepackage{setspace}
\usepackage{float}
\usetikzlibrary{calc}
\usepackage[super]{nth}
\usepackage{tikz}
\usepackage{extarrows}
\usepackage{pgfplots}
\pgfplotsset{compat=1.15}
\usetikzlibrary{arrows}

\theoremstyle{definition}

\theoremstyle{remark}

\numberwithin{equation}{section}

\allowdisplaybreaks

\newcommand\numberthis{\addtocounter{equation}{1}\tag{\theequation}}
\newcommand\restr[2]{{
  \left.\kern-\nulldelimiterspace 
  #1 
  \littletaller 
  \right|_{#2} 
  }}
\newcommand{\littletaller}{\mathchoice{\vphantom{\big|}}{}{}{}}  

\begin{document}

\title[]{SOME $Q$-HYPERGEOMETRIC IDENTITIES ASSOCIATED WITH PARTITION THEOREMS OF LEBESGUE, SCHUR AND CAPPARELLI}

\author{Yazan Alamoudi}
\address{Department of Mathematics, University of Florida, Gainesville,
FL 32611, USA}
\email{yazanalamoudi(at)ufl.edu}
\author{Krishnaswami Alladi}
\address{Department of Mathematics, University of Florida, Gainesville,
FL 32611, USA}
\email{alladik(at)ufl.edu}

\dedicatory{Dedicated to Mourad Ismail on the occasion of his \nth{80} birthday}

\date{\today}



\begin{abstract}
Here, we establish a polynomial identity in three variables $a, b, c$, and with the degree of the polynomial given in terms of two integers $L, M$. By letting $L$ and $M$ tend to infinity, we get the 1993 Alladi-Gordon $q$-hypergeometric key-identity for the generalized Schur Theorem as well as the fundamental Lebesgue identity by two different choices of the variables. This polynomial identity provides a generalization and a unified approach to the Schur and Lebesgue theorems. We discuss other analytic identities for the Lebesgue and Schur theorems and also provide a key identity ($q$-hypergeometric) for Andrews' deep refinement of the Alladi-Schur theorem. Finally, we discuss a new infinite hierarchy of identities, the first three of which relate to the partition theorems of Euler, Lebesgue, and Capparelli, and provide their polynomial versions as well.
\end{abstract}
\maketitle

\section{Introduction}\label{s1}
One of the fundamental $q$-hypergeometric identities is Lebesgue's identity:\\
\begin{align*}
\sum^{\infty}_{i=0}\frac{q^{T_i}(-cq)_i}{(q)_i}=\prod^{\infty}_{m=1}(1+q^m)(1+cq^{2m})=\prod^{\infty}_{m=1}\frac{(1+cq^{2m})}{(1-q^{2m-1})}.\numberthis\label{eq11}\\    
\end{align*}
The importance of \eqref{eq11} is due to the fact that when $c=0$ it yields Euler's series and product generating functions for partitions into distinct parts, and with the
dilations and translations given by\\
\begin{align*}
\text{dilation} \, q\mapsto q^2, \quad \text{translations} \, c\mapsto cq^{-1}
\, \text{or}\, c\mapsto cq,\numberthis\label{eq12}\\    
\end{align*}
it yields the $q$-hypergeometric identities for the Little G\"ollnitz partition theorems (see Theorem G below). In \cite{KABG93JCTA}, Alladi-Gordon gave the following combinatorial interpretation of Lebesgue's identity as a weighted partition theorem along with a combinatorial proof:\\

\textbf{Theorem L:}
\textit{Let $D(n;j)$ denote the number of partitions of $\pi: b_1+b_2+\cdots+b_{\nu}=n$ into distinct parts $b_i$, such that there are $j$ gaps $b_i-b_{i+1}\ge 2$ among the parts for $i=1,2,\cdots \nu$, with the convention $b_{\nu+1}=0$.}

\textit{Let $C(n;k)$ denote the number of partitions of $n$ with even parts non-repeating, such that there are precisely $k$ even parts. Then}\\
\begin{align*}
\sum_jD(n,j)(1+c)^j=\sum_kC(n;k)c^k.\\    
\end{align*}
Under the transformations in \eqref{eq12}, the two G\"ollnitz identities that emerge are:\\
\begin{align*}\label{eq12a}
\sum^{\infty}_{i=0}\frac{q^{i^2+i}(-cq^{-1};q^2)_i}{(q^2;q^2)_i}
=\prod^{\infty}_{m=1}(1+q^{4m})(1+q^{4m-2})(1+cq^{4m-3}),\tag{1.2a}    
\end{align*}\\
and
\begin{align*}\label{eq12b}
\sum^{\infty}_{i=0}\frac{q^{i^2+i}(-cq;q^2)_i}{(q^2;q^2)_i}
=\prod^{\infty}_{m=1}(1+q^{4m})(1+q^{4m-2})(1+cq^{4m-1}).\tag{1.2b}\\    
\end{align*}
The Little G\"ollnitz theorem(s) \cite{HG67}, which are the partition interpretations of \eqref{eq12a} and \eqref{eq12b}, are:\\

\textbf{Theorem G:}

\textit{For $i=1,2$, let $g_i(n;k)$ denote the number of partitions of $n$ into parts that differ by $\ge 2$, with strict inequality if a part is odd,  having $k$ odd parts, where the smallest part is $\ge i$.}

\textit{For $i=1,2$, let $G_i(n;k)$ denote the number of partitions into distinct parts, which are of the form $2,4,\, \text{or } 2i-1\pmod{4}$, and with $k$ parts $\equiv 2i-1\pmod{4}$. Then}\\
\begin{align*}
g_i(n;k)=G_i(n;k), \quad \text{for}\quad i=1,2.\\
\end{align*}

\textbf{Remark 1.1:}
It is only at the undilated level, that is, for Lebesgue's identity, the partition theorem (Theorem L) is a weighted partition theorem. Once we have the dilation as in \eqref{eq12}, Theorem G is a regular partition theorem (not weighted).\\

The celebrated 1926 partition theorem of Schur is:\\

\textbf{Theorem S:}

\textit{Let $B(n)$ denote the number of partitions of $n$ into distinct parts $\equiv\pm 1\pmod{3}$.}

\textit{Let $S(n)$ denote the number of partitions of $n$ into parts that differ by at least 3, with strict inequality if a part is a multiple of 3. Then}\\
\begin{align*}
S(n)=B(n).\\
\end{align*}

\textbf{Remark 1.2:}
Note one similarity in the difference conditions in Theorems G and S, namely, in Theorem G, the gap between parts is $\ge 2$ with strict inequality if a part is odd, while in Theorem S, the gap between parts is $\ge 3$ with strict inequality if a part is a multiple of $3$. For the remainder of this manuscript, we will refer to the partitions enumerated by $S(n)$ as \textit{Schur partitions}.\\

Gleissburg \cite{WG28} showed that Theorem S can be refined to\\
\begin{align*}
B(n;k)=S(n;k),\\
\end{align*}
where $B(n;k)$ and $S(n;k)$ denote the number of partitions enumerated by $B(n)$ and $S(n)$ with the condition that the number of parts is $k$, and with the convention that parts which are multiples of 3 are counted twice by $S(n;k)$.\\

In 1993, Alladi and Gordon \cite{KABG93MM} proved a two-parameter refinement and generalization of Theorem S, and in doing so, for the first time, cast Theorem S in the form of a $q$-hypergeometric identity, which they dubbed a \textit{key-identity}:\\
\begin{align*}
\sum_{\alpha, \beta, \gamma}\frac{a^{\alpha+\gamma}b^{\beta+\gamma}q^{T_{s}+T_{\gamma}}}{(q)_{\alpha}(q)_{\beta}(q)_{\gamma}}=\sum_{i,j}\frac{a^ib^jq^{T_i+T_j}}{(q)_i(q)_j}=\prod_{i=1}^{\infty}(1+aq^i)(1+bq^i),\numberthis\label{eq13}\\    
\end{align*}
where $s=\alpha+\beta+\gamma$.\\

By using the transformations\\
\begin{align*}
(\text{dilation}) \, q\mapsto q^3, \quad \text{and} \quad (\text{translations}) \, a\mapsto aq^{-2}, \, b\mapsto bq^{-1}\numberthis\label{eq14}\\
\end{align*}
in \eqref{eq13}, the following strong refinement of Theorem S falls out:\\

\textbf{Theorem A-G:}

\textit{Let $B(n;i,j)$ denote the number of partitions of $n$ into $i$ distinct parts $\equiv 1\pmod{3}$, and $j$ distinct parts $\equiv 2\pmod{3}$.}

\textit{Let $S(n;\alpha, \beta, \gamma)$ denote the number of partitions of the type enumerated by $S(n)$, with the condition that the number of parts $\equiv 1,2,3\pmod{3}$ is $\alpha, \beta, \gamma$ respectively. Then}\\
\begin{align*}
\sum_{\alpha+\gamma=i, \beta+\gamma=j}S(n;\alpha, \beta, \gamma)=B(n;i,j).\\
\end{align*}

Notice that in Theorem A-G, the total number of parts is\\
\begin{align*}
i+j=\alpha+\gamma+\beta+\gamma=\alpha+\beta+2\gamma,\\
\end{align*}
and so the parts that are multiples of $3$ are counted twice. In \cite{KABG93MM}, the combinatorial interpretation of the key-identity was given in terms of partitions into parts occurring in three possible colors: two primary colors, $a$ and $b$, and one secondary color, $ab$ the combination of the other two, with gap conditions on the colored parts, and with the condition that the parts occurring in the secondary color are counted twice. Under the dilation and translations given in (1.4), the parts in primary colors $a,b$ correspond to parts $\equiv 1,2\pmod{3}$, and parts in secondary color are then the multiples of 3.\\

The colored partition version (generalization) of Theorem A-G is proved combinatorially (bijectively) in \cite{KABG93MM}. The combinatorial proofs of the weighted partition Theorem L given in \cite{KABG93JCTA}, and of the colored partition version of Theorem A-G in \cite{KABG93MM} are similar, with the main difference being in the final step, where in the case of Theorem L, a certain choice could be made; this is why $D(n;j)$ has a weight $(1+c)^j$ attached to it. Since the combinatorial proofs of Theorem L and the colored generalization of Theorem A-G are so similar, it is natural to ask if there is a unified $q$-hypergeometric approach to Lebesgue's identity \eqref{eq11} and the key-identity \eqref{eq13} for the generalized Schur theorem? After establishing a finite analog of Lebesgue's identity in Section \ref{s3}, we provide in Section \ref{s6} a new polynomial identity, from which, under two different specializations, the finite Schur and Lebesgue identities fall out. Following this, in Section \ref{s6}, we provide a $q$-hypergeometric \textit{key-identity} for Andrews' deep refinement of the Alladi-Schur Theorem. Finally, in Section \ref{s7}, we introduce a new infinite hierarchy of identities of which the first three correspond to the partition theorems of Euler, Lebesgue, and Capparelli; we provide a polynomial version of this infinite hierarchy as well.\\

To aid the reader, we will conclude this introductory section by recalling basic facts, along with notational conventions, that are used throughout the paper.\\

For complex numbers $a, q$, we use the $q$-Pochhammer symbols
\begin{align*}
(a)_n=(a;q)_n:=\prod^{n-1}_{j=0}(1-aq^j),
\end{align*}
and
\begin{align*}
(a;q)_{\infty}=\lim_{n\to\infty}(a;q)_j=\prod^{\infty}_{j=0}(1-aq^j),\quad \text{if} \quad
|q|<1.
\end{align*}
The variable $q$ is called the base. We often write $(a)_n$ in place of
$(a;q)_n$ suppressing $q$, but when the base is anything other than $q$,
it will be displayed.\\

We also make use of the $q$-binomial coefficients given by\\
\begin{align*}
{n\brack m}= {n\brack m}_q:=\frac{(q)_n}{(q)_m(q)_{n-m}}, \quad \text{for} \quad 0\le m\le n\numberthis\label{eq15}\\
\end{align*}
which are polynomials in $q$ of degree $m(n-m)$. When $n\geq0$, the $q$-binomial coefficients have value 0 when $m<0$ or when $m>n$. This is because $1/(q)_j=0$ when $j<0$.\\

We shall often use the following identity involving the $q$-binomial coefficients:\\
\begin{align*}
(-cq)_n=\sum_{k\geq0}c^kq^{T_k}{n\brack k},\numberthis\label{eq16}\\
\end{align*}
where, here and throughout, $T_k=k(k+1)/2$ is the $k$-th Triangular number. In some instances, the expressions involving the $q$-Pochhammer symbol and the $q$-binomial coefficients presented in this paper remain valid with $n<0$. In such cases, the meaning of ${n\brack k}$ and $(-cq)_n$ is as \cite{SFAS19}. Furthermore, \eqref{eq16} remains valid for any integer $n$ (even $n<0$) as can be seen from \cite[Thm 4.10]{SFAS19}. \\

\section{A very short proof of the key identity for Schur's theorem}\label{s2}
The proof of the key identity \eqref{eq13} given in \cite{KABG93MM} utilizes the $q$-Chu-Vandermonde summation. A second proof of \eqref{eq13} was given in \cite{KABG93MM} by rewriting it suitably and then using a Durfee rectangles argument. In October 2022, the second author communicated \cite{KA22} to George Andrews a very short proof of the key identity \eqref{eq13}, which we give here because this proof and the underlying combinatorics motivate the construction of the polynomial identity that provides the unification:\\

\textbf{Proof:} Begin by expanding $(-aq)_{\infty}$ and by splitting $(-bq)_{\infty}$ as follows:\\
\begin{align*}
(-aq)_{\infty}(-bq)_{\infty}=\sum^{\infty}_{i=0}\frac{a^iq^{T_i}}{(q)_i}(-bq)_i(-bq^{i+1})_{\infty}.\numberthis\label{eq21}\\
\end{align*}
Next, expand $(-bq)_i$ and $(-bq^{i+1})_{\infty}$, and substitute these expansions in \eqref{eq21} to get\\
\begin{align*}
(-aq)_{\infty}(-bq)_{\infty}&=\sum^{\infty}_{i=0}\frac{a^iq^{T_i}}{(q)_i}\left(\sum^i_{j=0}b^jq^{T_j}{i\brack j}\right)\left(\sum^{\infty}_{\ell=0}\frac{b^{\ell}q^{T_{\ell}+i\ell}}{(q)_{\ell}}\right)\\
&=\sum^{\infty}_{i=0}a^iq^{T_i}\left(\sum^i_{j=0}\frac{b^jq^{T_j}}{(q)_j(q)_{i-j}}\right)\left(\sum^{\infty}_{\ell=0}\frac{b^{\ell}q^{T_{\ell}+i\ell}}{(q)_{\ell}}\right).\numberthis\label{eq22}\\
\end{align*}
At this stage, consider the following replacements in \eqref{eq22}\\
\begin{align*}
j\mapsto\gamma, \quad i-j=i-\gamma\mapsto\alpha, \quad \text{and} \quad
\ell\mapsto \beta,\numberthis\label{eq23}\\
\end{align*}
to rewrite the expression on the right in \eqref{eq22} as\\
\begin{align*}
\sum_{\alpha,\beta,\gamma}\frac{a^{\alpha+\gamma}b^{\beta+\gamma}q^{T_{s}+T_{\gamma}}}{(q)_{\alpha}(q)_{\beta}(q)_{\gamma}},\numberthis\label{eq24}\\    
\end{align*}
where we have used the identity\\
\begin{align*}
T_n+T_m+nm=T_{n+m}\\  
\end{align*}
for Triangular numbers, and $s=\alpha+\beta+\gamma$. The key identity follows from \eqref{eq21} and \eqref{eq24}.\qed\\

\textbf{Remark 2.1:}
In \cite{KABG93MM}, the first two steps in the combinatorial proof of the Generalized Schur Theorem were as follows: Start with a vector partition $<\pi_a, \pi_b>$ in which $\pi_a$ is a partition into $i$ parts in color $a$, all distinct, and $\pi_b$ is a partition into $j$ parts in color $b$, all distinct. Then separate the parts of $\pi_b$ into those that are $\le i$ in size and those that are $>(i+1)$ in size. There are six steps in that combinatorial proof, but these first two steps correspond to\\
\begin{align*}
\frac{a^iq^{T_i}}{(q)_i}(-bq)_i(-bq^{i+1})_{\infty},\\
\end{align*}
and this motivated the starting point of the short proof of the key identity.\\

\section{A finite version of Lebesgue's identity}\label{s3}
The $q$-binomial coefficients ${n\brack m}$ have the property that\\
\begin{align*}
\lim_{n\to\infty}{n\brack m}=\frac{1}{(q)_m}.\numberthis\label{eq31}\\
\end{align*}
So, a natural way to construct polynomial analogs of $q$-hypergeometric identities is to bring in $q$-binomial coefficients in place of terms like $1/(q)_m$. We now establish a polynomial version of Lebesgue's identity:\footnote{An equivalent identity appears in \cite{MJR10} in a different form. However, our approach is different.}\\

\textbf{Theorem 3.1:}
\textit{For all positive integers $M$, we have}\\
\begin{align*}
\sum^M_{i=0}q^{T_i}(-cq)_i{M\brack i}=(-q)_M\left(\sum^M_{k=0}\frac{c^kq^{2T_k}}{(q^2;q^2)_k}(q^{M-k+1})_k\right).\\
\end{align*}

\textbf{Proof:}
Begin by expanding $(-cq)_i$ and using \eqref{eq15} to rewrite the left-hand side of the expression in Theorem 3.1 as\\
\begin{align*}
\sum^M_{i=0}q^{T_i}(-cq)_i{M\brack i}=\sum^M_{i=0}q^{T_i}\left(\sum^i_{k=0}c^kq^{T_k}{i\brack k}\right){M\brack i}=\sum_{i,k}c^kq^{T_i+T_k}\frac{(q)_M}{(q)_k(q)_{i-k}(q)_{M-i}}.\numberthis\label{eq32}\\    
\end{align*}
If we set $I=i-k$, we may rewrite the right-hand side of \eqref{eq32} as\\
\begin{align*}
\sum_{I,k}c^kq^{T_{I+k}+T_k}\frac{(q)_M}{(q)_k(q)_I(q)_{M-I-k}}=\sum_{I,k}c^kq^{T_{I+k}+T_k}{{M-I}\brack k}{M\brack I}\numberthis\label{eq33}\\ 
\end{align*}
again by \eqref{eq15}. At this stage, we replace $T_{I+k}$ in \eqref{eq33} with $T_I+T_k+Ik$ to rewrite \eqref{eq33} as\\
\begin{align*}
\sum_{I,k}c^kq^{2T_k+T_I+Ik}{{M-I}\brack k}{M\brack I}=\sum_{I,k}c^kq^{2T_k+T_I+Ik}{M\brack k}{{M-k}\brack I},\numberthis\label{eq34}\\
\end{align*}
using \eqref{eq15} once more. Finally, we write the right-hand side of \eqref{eq34} as\\
\begin{align*}
&\sum_kc^kq^{2T_k}{M\brack k}\left(\sum_Iq^{T_I+Ik}{{M-k}\brack I}\right)=\sum_kc^kq^{2T_k}{M\brack k}(-q^{k+1})_{M-k}\\
&=(-q)_M\sum^M_{k=0}\frac{c^kq^{2T_k}}{(q^2;q^2)_k}\frac{(q)_M}{(q)_{M-k}}=(-q)_M\sum^M_{k=0}\frac{c^kq^{2T_k}}{(q^2;q^2)_k}(q^{M-k+1})_k,\numberthis\label{eq35}\\
\end{align*}
and this proves Theorem 3.1.\qed\\ 


{\bf{Remark 3.2:}} Ole Warnaar has pointed out (private correspondence, 2025) that
Theorem 3.1 is a special case of one of Jackson's classic ${}_2\phi_1$ transformations, by setting $c=0, a=q^{-n}$ and replacing $z$ by $zq^n$ in Gasper and Rahman
\cite[III.4]{GGMR90}. Our emphasis has been to provide direct proofs of this and other identities in this paper.\medskip

\newpage

\textbf{Lebesgue's identity as a limiting case of Theorem 3.1:}\\

Let $M\to\infty$ in Theorem 3.1. Then, in view of \eqref{eq31}, the left-hand side of Theorem 3.1 is\\
\begin{align*}
\sum^{\infty}_{i=0}\frac{q^{T_i}(-cq)_i}{(q)_i},\numberthis\label{eq36}\\
\end{align*}
which is the left-hand side of \eqref{eq11}. On the other hand, when $M\to\infty$, the right-hand side of Theorem 3.1 becomes\\
\begin{align*}
(-q)_{\infty}\sum^{\infty}_{k=0}\frac{c^kq^{2T_k}}{(q^2;q^2)_k}=(-q)_{\infty}(-cq^2;q^2)_{\infty},\numberthis\label{eq37}\\
\end{align*}
because $(q^{M-k+1})_k\to 1$ as $M\to\infty$. This yields Lebesgue's identity.\\

\textbf{Another finite version of Lebesgue's identity:}\\\quad\\

There are several possible finite versions of Lebesgue's identity, such as\\
\begin{align*}
\sum^{m+n}_{N=0}q^{(N^2+N)/2}\sum^N_{k=0}b^kq^{(k^2+k)/2}{m\brack k}{n\brack {N-k}}=(-q)_n(q)_m\sum^m_{i=0}\frac{b^iq^{i^2+i}(-q^{n+1})_i}{(q^2;q^2)_i(q)_{m-i}}
\numberthis\label{eq38}\\
\end{align*}
which is due to Alladi (1994 - unpublished), but we have emphasized the finite version in Theorem 3.1 because this is connected to the unified approach to the Schur and Lebesgue identities that will be given below.\\

To realize that \eqref{eq38} is a finite version of Lebesgue's identity, let $m,n \to\infty$ in \eqref{eq38}. Then the right-hand side \eqref{eq38} becomes\\
\begin{align*}
(-q)_{\infty}(q)_{\infty}\sum^{\infty}_{i=0}\frac{b^iq^{i^2+i}}{(q^2;q^2)_i(q)_{\infty}}=(-q)_{\infty}(-bq^2;q^2)_{\infty},\numberthis\label{eq39}\\
\end{align*}
which is the right-hand side of \eqref{eq11}. Under these limits, the left side of \eqref{eq38} becomes\\
\begin{align*}
\sum^{\infty}_{N=0}q^{T_N}\sum^N_{k=0}\frac{b^kq^{T_k}}{(q)_k(q)_{N-k}}&=\sum^{\infty}_{N=0}\frac{q^{T_N}}{(q)_N}\sum^N_{k=0}b^kq^{T_k}{N\brack k}=\sum^{\infty}_{N=0}\frac{q^{T_N}(-bq)_N}{(q)_N},\numberthis\label{eq310}\\
\end{align*}
which is the left-hand side of \eqref{eq11}.\\

We now sketch the proof of \eqref{eq38} which is non-trivial.\\

\textbf{Proof of (3.8):}
In Alladi \cite{KA94}, the following was established both hypergeometrically and combinatorially:\\

\textbf{Lemma 3.3:} (Transformation Formula)\\
\begin{align*}
\sum^{\infty}_{n=0}\frac{a^nq^{n^2}(-bq;q^2)_n}{(q^2;q^2)_n}=\sum^{\infty}_{n=0}\frac{(ab)^nq^{2n^2}(-aq^{2n+1};q^2)_{\infty}}{(q^2;q^2)_n}.\\
\end{align*}

To prove the Lemma, expand $(-bq;q^2)_n$ on the left-hand side using \eqref{eq16} and reverse the order of summation to get the right-hand side. The combinatorial proof of Lemma 3.3 uses a redistribution idea of Bressoud (see \cite{DMB78,DMB79}).\\

The finite version of Lemma 3.3 is\\

\textbf{Lemma 3.4:}\\
\begin{align*}
\sum^{m+n}_{N=0}a^Nq^{N^2}\sum^n_{k=0}b^kq^{k^2}{m\brack k}_{q^2}{n\brack {N-k}}_{q^2}=\sum^m_{i=0}(ab)^iq^{2i^2}{m\brack i}_{q^2}(-aq^{2i+1};q^2)_n.\\
\end{align*}

To prove Lemma 3.4, expand $(-aq^{2i+1};q^2)_n$ to rewrite the right-hand side as\\
\begin{align*}
\sum^m_{i=0}(ab)^iq^{2i^2}{m\brack i}_{q^2}\sum^n_{j=0}a^jq^{j^2+2ij}{n\brack j}_{q^2}.\numberthis\label{eq311}\\
\end{align*}
If we rearrange the right-hand side of \eqref{eq311} by putting $i+j=N$, we get the left-hand side of Lemma 3.4, thereby proving it.\\

Next in Lemma 3.4, replace $a$ by $aq$ and $b$ by $bq$ to get\\
\begin{align*}
\sum^{m+n}_{N=0}a^Nq^{N^2+N}\sum^N_{k=0}b^kq^{k^2+k}{m\brack k}_{q^2}{N\brack {N-k}}_{q^2}=\sum^m_{i=0}(ab)^iq^{2i^2+2i}{m\brack i}_{q^2}(-aq^{2i+2};q^2)_n.\numberthis\label{eq312}\\   
\end{align*}
In (3.12) replace $q^2$ with $q$ to get\\
\begin{align*}
\sum^{m+n}_{N=0}a^Nq^{T_N}\sum^n_{k=0}b^kq^{T_k}{m\brack k}{n\brack {N-k}}=\sum^m_{i=0}(ab)^iq^{i^2+i}{m\brack i}(-aq^{i+1})_n.\numberthis\label{eq313}\\    
\end{align*}
If we set $a=1$ in \eqref{eq313}, the right-hand side becomes\\
\begin{align*}
(q)_m\sum^m_{i=0}\frac{b^iq^{i^2+i}(-q^{i+1})_n}{(q)_i(q)_{m-i}}=(-q)_n(q)_m\sum^m_{i=0}\frac{b^iq^{i^2+i}(-q^{n+1})_i}{(q^2;q^2)_i(q)_{m-i}},\numberthis\label{eq314}\\    
\end{align*}
and so \eqref{eq313} and \eqref{eq314} yield \eqref{eq38}.\qed\\

\section{A power series identity unifying Schur and Lebesgue}\label{s4}
Here we prove the following power series identity due to Alamoudi, from which \eqref{eq11} and \eqref{eq13} emerge as special cases:\\

\textbf{Theorem 4.1:}
\textit{With free parameters $a,b,c$, we have}\\
\begin{align*}
\sum^{\infty}_{i=0}\frac{a^iq^{T_i}}{(q)_i}\left(-\frac{c}{a}q\right)_i(-bq^{i+1})_{\infty}=\sum_{t,j,\ell}\frac{a^tb^{\ell}c^jq^{T_{t+\ell+j}+T_j}}{(q)_t(q)_{\ell}(q)_j}.\\
\end{align*}

\textbf{Proof:}
We expand $(-\frac{c}{a}q)_i$ and $(-bq^{i+1})_{\infty}$ to get\\
\begin{align*}
\sum^{\infty}_{i=0}\frac{a^iq^{T_i}}{(q)_i}\left(-\frac{c}{a}q\right)_i(-bq^{i+1})_{\infty}&=\sum^{\infty}_{i=0}\frac{a^iq^{T_i}}{(q)_i}\left(\sum^i_{j=0}\left(\frac{c}{a}\right)^jq^{T_j}{i\brack j}\right)\left(\sum^{\infty}_{\ell=0}\frac{b^{\ell}q^{T_{\ell}+i{\ell}}}{(q)_{\ell}}\right)\\
&=\sum^{\infty}_{i=0}a^iq^{T_i}\left(\sum^i_{j=0}\frac{(\frac{c}{a})^jq^{T_j}}{(q)_j(q)_{i-j}}\right)\left(\sum^{\infty}_{\ell=0}\frac{b^{\ell}q^{T_{\ell}+i{\ell}}}{(q)_{\ell}}\right).\numberthis\label{eq41}\\
\end{align*}
Now, for the sum on the right in \eqref{eq41}, put $t=i-j$ and simplify to get\\
\begin{align*}
\sum_{t,j,\ell}\frac{a^tb^{\ell}c^{j}q^{T_{t+j}+T_j+T_{\ell}+(t+j){\ell}}}{(q)_t(q)_j(q)_{\ell}}=\sum_{t,j,\ell}\frac{a^tb^{\ell}c^{j}q^{T_{t+j+\ell}+T_j}}{(q)_t(q)_j(q)_{\ell}},\\   
\end{align*}
which proves Theorem 4.1.\\ 

We record two corollaries to Theorem 4.1:\\

\textbf{Corollary 4.2:}
\textit{The key-identity \eqref{eq13} for the generalized Schur's theorem holds}.\\

\textbf{Proof:}
In Theorem 4.1, take $c=ab$. Then, the left-hand side of Theorem 4.1 is\\
\begin{align*}
(-bq)_{\infty}\sum^{\infty}_{i=0}\frac{a^iq^{T_i}}{(q)_i}=(-aq)_{\infty}(-bq)_{\infty}.\numberthis\label{eq42}\\
\end{align*}
The right-hand side of Theorem 4.1 is\\
\begin{align*}
\sum_{t,j,\ell}\frac{a^{t+j}b^{\ell+j}q^{T_{t+j+\ell}+T_j}}{(q)_t(q)_j(q)_{\ell}}.\numberthis\label{eq43}\\
\end{align*}
Now, \eqref{eq13} follows from \eqref{eq42} and \eqref{eq43} with the replacements
\begin{align*}\\
t\mapsto\alpha, \quad \ell\mapsto\beta, \quad \text{and} \quad j\mapsto\gamma.
\end{align*}\\
Hence Corollary 4.2.\qed\\

\textbf{Corollary 4.3:}
\textit{The Lebesgue identity \eqref{eq11} holds}.\\

\textbf{Proof:}
Take $b=0, a=1$ in Theorem 4.1. Then, the left-hand side of Theorem 4.1 is the left-hand side of \eqref{eq11}. Since $b=0$, the only contribution to the right-hand side of Theorem 4.1 is from $\ell=0$, interpreting $b^0=1$ always. So the right-hand side replacing $t\mapsto i$ is\\
\begin{align*}
\sum_{i.j}\frac{c^jq^{T_{i+j}+T_j}}{(q)_i(q)_j}&=\sum_{i,j}\frac{c^jq^{T_i+2T_j+ij}}{(q)_i(q)_j}\\
&=\sum^{\infty}_{j=0}\frac{c^jq^{2T_j}}{(q)_j}\left(\sum^{\infty}_{i=0}\frac{q^{T_i+ij}}{(q)_j}\right)=\sum^{\infty}_{j=0}\frac{c^jq^{2T_j}}{(q)_j}(-q^{j+1})_{\infty}\\
&=(-q)_{\infty}\sum^{\infty}_{j=0}\frac{c^jq^{2T_j}}{(q^2;q^2)_j}=(-q)_{\infty}(-cq^2;q^2)_{\infty},\\
\end{align*}
and this yields \eqref{eq11}. Hence Corollary 4.3.\qed\\

\textbf{Remark 4.4:}
Since the Schur key-identity and the Lebesgue identity fall out as corollaries (but as two different special cases), Theorem 4.1 provides the unification stressed at the beginning of the paper. This unification is facilitated by the introduction of a third \textit{free parameter} $c$ in Theorem 4.1. In \cite{KABG93MM} where the key-identity is proved, the symbol $c$ is used to denote parts of secondary color, but $c$ was always taken to be $ab$ to get the product on the right in \eqref{eq13}. The motivation to introduce the third parameter $c$ in Theorem 4.1 is from the short proof of the key identity; the split product
\begin{align*}
(-bq)_i(-bq^{i+1})_{\infty}\\
\end{align*}
in the short proof is replaced in Theorem 4.1 by the more general split product\\
\begin{align*}
\left(-\frac{c}{a}q\right)_i(-bq^{i+1})_{\infty}.\\
\end{align*}
The above coincide when $c=ab$. Now, again, consider the bijection in \cite{KABG93MM}. The sub-partition $\lambda_b$ of $\pi_b$ containing the parts $\leq\nu(\pi_a)$ in size is used to construct a new partition $\pi_{a,c}=\pi_a+(\lambda_b)^*$ with $\nu(\pi_{a,c})=\nu(\pi_a)$ but now $\nu(\lambda_b)$ of the parts have become $c$ parts.\footnote{$\pi^*$ denotes the conjugate of $\pi$.} This motivates the $\frac{c}{a}$ factor.\\

\textbf{Remark 4.5:}
Ramamani and Venkatachaliengar \cite{VRKV72} generalized Lebesgue's identity \eqref{eq11} as follows:\\
\begin{align*}
\sum^{\infty}_{i=0}\frac{t^iq^{T_i}(z)_i}{(q)_i}=(z)_{\infty}(-tq)_{\infty}\sum^{\infty}_{j=0}\frac{z^j}{(q)_j(-tq)_j}.\numberthis\label{eq44}\\
\end{align*}
Identity \eqref{eq44} can be proved $q$-hypergeometrically or combinatorially using vector partitions. Now \eqref{eq11} can be deduced from \eqref{eq44} as follows: Take $t=1$ and $z=-cq$. Then the left-hand side of \eqref{eq44} is the left-hand side of \eqref{eq11}. With these values of $z$ and $t$, the right-hand side of \eqref{eq44} is\\
\begin{align*}
(-cq)_{\infty}(-q)_{\infty}\sum^{\infty}_{j=0}\frac{(-cq)^j}{(q^2;q^2)_j}=\frac{(-cq)_{\infty}(-q)_{\infty}}{(-cq;q^2)_{\infty}}=(-q)_{\infty}(-cq^2;q^2)_{\infty},\\
\end{align*}
which is the right-hand side of \eqref{eq11}. Thus Lebesgue's identity follows from \eqref{eq44}, but is different from our derivation of \eqref{eq11} from Theorem 4.1, because we get $(-cq^2;q^2)_{\infty}$ directly, whereas from \eqref{eq44}, $(-cq^2;q^2)_{\infty}$ is obtained from the cancellation in\\
\begin{align*}
\frac{(-cq)_{\infty}}{(-cq;q^2)_{\infty}}.\\
\end{align*}

\textbf{Remark 4.6:}
The replacement $c\mapsto cq^{-1}$ in \eqref{eq11} yields the equivalent identity\\
\begin{align*}
\sum^{\infty}_{i=0}\frac{q^{T_i}(-c)_i}{(q)_i}=\prod^{\infty}_{m=1}(1+q^m)(1+cq^{2m-1})=\prod^{\infty}_{m=1}\frac{(1+cq^{2m-1})}{(1-q^{2m-1})}.\numberthis\label{eq45}\\  
\end{align*}
whose combinatorial interpretation yields Sylvester's famous refinement \cite{JJS82} of Euler's theorem. Ramamani and Venkatachaliengar actually generalize \eqref{eq45} by establishing an identity equivalent to \eqref{eq44}; we have preferred the version of their identity as in \eqref{eq44} in view of the discussion of Lebesgue's identity in this paper.\\

\textbf{Remark 4.7:} Since the right-hand side in Theorem 4.1 is symmetric under interchanging\\ $a$ and $b$, it follows that 

\begin{align*}
\sum^{\infty}_{i=0}\frac{a^iq^{T_i}}{(q)_i}\left(-\frac{c}{a}q\right)_i(-bq^{i+1})_{\infty}=\sum^{\infty}_{i=0}\frac{b^iq^{T_i}}{(q)_i}\left(-\frac{c}{b}q\right)_i(-aq^{i+1})_{\infty}.\\
\end{align*}

\textbf{Remark 4.8:}
We alert the reader that the dilation and translations in \eqref{eq14}, as well as the translation $c\to cq^{-3}$, the right-hand side of Theorem 4.1 is the sum over generating functions of Schur partitions with the powers of $a$ counting the parts $\equiv 1\pmod{3}$, the powers of $b$ counting the parts $\equiv 2\pmod{3}$ and the powers of $c$ counting the parts divisible by $3$. More generally, the function
\begin{align*}
G(t,\ell,j)=\frac{q^{T_{t+\ell+j}+T_j}}{(q)_t(q)_{\ell}(q)_j}
\end{align*}
counts \textit{Type-1} partitions of a prescribed number of parts of each color, specifically, $t$ $a$-parts, $\ell$ $b$-parts, and $j$ $c$-parts (see \cite{KABG95}).  \textit{Type-1} refers to a general class of colored partitions, whose exact definition is given in \cite{KABG95}, that amount to Schur partitions under standard transformations. Furthermore, in  \cite{KABG95}, Alladi-Gordon demonstrated that there are six schemes (i.e., Type-2 up to Type-6), all counted by $G(t,\ell,j)$.\\

In the next section, we shall establish a polynomial version of Theorem 4.1. 

\section{A general polynomial identity in three parameters}\label{s5}
In this section, we prove a general multi-parameter polynomial identity due to Alamoudi, from which some of the key results stated above follow either as limiting cases or as special cases.\\

\textbf{Theorem 5.1:} (\textit{Finite three-parameter Schur})
\textit{For any pair of integers $L, M$, and parameters $a,b,c$, we have}\\
\begin{align*}
\sum_{i\geq0}a^iq^{T_i}\left(-\frac{c}{a}q\right)_i(-bq^{i+1})_{L-i}{M\brack i}=\sum_{i,j,k\geq0}a^ib^jc^kq^{T_{i+j+k}+T_k}{{M-i}\brack k}{M\brack i}{{L-i-k}\brack j}.\\
\end{align*}

{\bf{Proof:}} Expand $\left(-\frac{c}{a}q\right)_i$ and $(-bq^{i+1})_{L-i}$ to rewrite the
left side of Theorem 5.1 as\\
\begin{align*}
&\sum_{i\geq0}a^iq^{T_i}\left(-\frac{c}{a}q\right)_i(-bq^{i+1})_{L-i}{M\brack i}\\
&=\sum_{i\geq0}a^iq^{T_i}\left(\sum^i_{k=0}(\frac{c}{a})q^{T_k}{i\brack k}\right)\left(\sum_{j\geq 0}b^jq^{T_j+ij}{{L-i}\brack j}\right){M\brack i}\\
&=\sum_{\substack{i,j,k\geq 0\\k\leq i }}a^{i-k}b^jc^kq^{T_{i+j}+T_k}{{L-i}\brack j}{i\brack k}{M\brack i}\\
&=\sum_{\substack{i,j,k\geq 0\\k\leq i }}a^{i-k}b^jc^kq^{T_{i+j}+T_k}{{L-i}\brack j}{{M-(i-k)}\brack k}{M\brack {i-k}},\numberthis\label{eq51}\\
\end{align*}
because\\
\begin{align*}
{i\brack k}{M\brack i}={{M-i+k}\brack k}{M\brack {i-k}},\numberthis\label{eq52}\\    
\end{align*}
Theorem 5.1 follows by replacing $i$ by $i+k$ in \eqref{eq51}.\qed\\

\textbf{Remark 5.2:}
Another way to write the right-hand side of Theorem 5.1 to make it more appealing combinatorially and symmetric is to replace\\
\begin{align*}
M\mapsto M'-j-k,\numberthis\label{eq53}\\
\end{align*}
which converts it to\\
\begin{align*}
\sum_{i,j,k}a^ib^jc^kq^{T_{i+j+k}+T_k}{{M'-(i+j+k)}\brack k}{{M'-(j+k)}\brack i}{{L-(i+k)}\brack j}.\numberthis\label{eq54}\\  
\end{align*}
We now consider the consequences of Theorem 5.1.\\  

\textbf{Corollary 5.3:}
\textit{Theorem 4.1 holds.}\\

\textbf{Proof:} Let $M,L\to\infty$ in Theorem 5.1 to get Theorem 4.1.\qed\\

\textbf{Corollary 5.4:}
\textit{The following finite (polynomial) version of Lebesgue's identity holds:}\\
\begin{align*}
\sum^{\infty}_{i=0}q^{T_i}(-cq)_i{M\brack i}=\sum_{i,k}c^kq^{T_i+2T_k+ik}{{M-i}\brack k}{M\brack i}\\
\end{align*}

\textbf{Proof:} Take $b=0$, $a=1$ in Theorem 5.1. Since $b=0$, the values $j>0$ do not make a contribution. Thus, we put $j=0$. This yields Corollary 5.4.\qed\\

\textbf{Remark 5.5:} The left-hand side of Corollary 5.4 is identical to the left-hand side
of Theorem 3.1. But the right-hand side of Corollary 5.4 is very different from the right-hand side of Theorem 3.1. Thus, Corollary 5.4 provides a different finite version of Lebesgue's identity. However, the right-hand side of Corollary 5.4 can be transformed into the right-hand side of Theorem 3.1. In fact, the left side of \eqref{eq34} is identical to the right-hand side of Corollary 5.4. \\

\textbf{Corollary 5.6:} \textit{Theorem 3.1 follows from Theorem 5.1.}\\

To realize that Corollary 5.4 is indeed a finite version of Lebesgue's identity \eqref{eq11}, let $M\to\infty$ in Corollary 5.4. Then, the left-hand side of Corollary 5.4 is clearly\\
\begin{align*}
\sum^{\infty}_{i=0}\frac{q^{T_i}(-cq)_i}{(q)_i}\\
\end{align*}
which is the left-hand side of (1.1). When $M\to\infty$, the right-hand side of
Corollary 5.4 becomes\\
\begin{align*}
&\sum_{i,k}\frac{c^kq^{T_i+2T_k+ik}}{(q)_i(q)_k}=\sum^{\infty}_{k=1}\frac{c^kq^{2T_k}}{(q)_k}\sum^{\infty}_{i=0}\frac{q^{T_i+ik}}{(q)_i}\\
&=\sum^{\infty}_{k=1}\frac{c^kq^{2T_k}}{(q)_k}(-q^{k+1})_{\infty}=(-q)_{\infty}\sum^{\infty}_{k=1}\frac{c^kq^{2T_k}}{(q^2;q^2)_k}=(-q)_{\infty}(-cq^2;q^2)_{\infty},\\
\end{align*}
which is the right-hand side of \eqref{eq11}.\\

\textbf{Remark 5.7:} We point out that the product of the two $q$-binomial coefficients in Corollary 5.4 can be rewritten as\\
\begin{align*}
{{M-i}\brack k}{M\brack i}={M\brack i,k,M-i-k}=\frac{(q)_M}{(q)_i(q)_k(q)_{M-i-k}},\numberthis\label{eq55}\\
\end{align*}
a $q$-multinomial\footnote{In the case of the $q$-multinomial coefficent of order $3$, such as in \eqref{eq55}, we display all three indices $i,k$ and $M-i-k$ whereas for the $q$-binomial coefficient ${M\brack i}$ we suppress $M-i$ (this is standard notation; see, for example, \cite{GEA76}).} coefficient of order 3. In \cite{KABG95}, Alladi-Gordon discuss how the generalized Schur partitions are related to $q$-multinomial coefficients of order 3.
This link of the Lebesgue identity with $q$-multinomial coefficients of order 3 is yet another Schur-Lebesgue unification. In Alladi-Berkovich [6; Eqn. 1.15], a finite version of Lebesgue's identity is established; that identity has the product of two $q$-binomial coefficients as in (5.5), but the link between Lebesgue's identity with $q$-multinomial coefficients of order $3$ and Schur's theorem is not considered in [6]. Warnaar \cite{SOW04} has provided a new proof of the Alladi-Berkovich finite version of the Lebesgue identity.\\


\textbf{Remark 5.8:} In \cite{KAAB02}, Alladi-Berkovich prove both combinatorially and $q$-theoretically the following double bounded version of the Alladi-Gordon key-identity  for Schur's partition theorem:\\
\begin{align*}
q^{T_i+T_j}{L\brack j}{{M-j}\brack i}=\sum_kq^{T_{i+j-k}+T_k}{{M-i-j+k}\brack k}{{M-j}\brack {i-k}}{{L-i}\brack {j-k}}
.\numberthis\label{eq56}\\
\end{align*}
If we multiply both sides of \eqref{eq56} by $a^ib^j$ and sum over $i,j$, we get a double bounded version of \eqref{eq13}. Our Theorem 5.1 has a third parameter $c$ and is different in shape from \eqref{eq56}, but with certain special choices and substitutions, Theorem 5.1 will yield \eqref{eq56}. In particular, if we set $c=ab$ in Theorem 5.1 and replace $i\to i'$ and $j\to j'$, then we get\\
\begin{align*}
(-aq)_M(-bq)_L 
    =\sum_{i',j',k}a^{i'}b^{j'}(ab)^kq^{T_{i'+j'+k}+T_k}{{M-i'}\brack k}{M\brack i'}{{L-i'-k}\brack j'}.\numberthis\label{eq57}\\
\end{align*}
In \eqref{eq57}, by setting $i=i'+k$, $j=j'+k$, and comparing the coefficients of $a^ib^j$ we get
\begin{align*}
\sum_{\substack{i',j',k\\i'+k=i,j'+k=j}}q^{T_{i'+j'+k}+T_k}{{M-i'}\brack k}{M\brack i'}{{L-i'-k}\brack j'}
=q^{T_i+T_j}{L\brack j}{{M}\brack i}.\numberthis\label{eq58}\\
\end{align*}
Setting $M\to M-j$ in \eqref{eq58} and simplifying the left-hand side yields \eqref{eq56}. Notice, however, that the double bounded version of \eqref{eq13} obtained by multiplying both sides of \eqref{eq56} by $a^ib^j$ and summing over $i,j$ is different from \eqref{eq57}.\\

\textbf{Remark 5.9:}
For $M\geq0$, setting $b=0$ in Theorem 5.1 gives
\begin{align*}
\sum_{i=0}^M a^iq^{T_i}(-\frac{c}{a}q)_i{ M\brack i}=\sum_{i,k}a^ic^k q^{T_{i+k}+T_k}{ M\brack i,k,M-i-k}.
\end{align*}
Letting $q\to 1$ gives 
\begin{align*}
\sum_{i=0}^M a^i(1+\frac{c}{a})^i{ M\choose i}=\sum_{i,k}a^ic^k { M\choose i,k,M-i-k}.
\end{align*}
The left-hand side is
\begin{align*}
\sum_{i=0}^M (a+c)^i{ M\choose i}=(1+a+c)^M.
\end{align*}
The above is an instance of the trinomial theorem.

\section{A key identity for the Alladi-Schur Theorem}\label{s6}
Schur's partition theorem has always been associated with the modulus $3$ or the
modulus $6$, the latter because\\
\begin{align*}
\prod^{\infty}_{m=1}(1+q^{3m-2})(1+q^{3m-1})=\prod^{\infty}_{m=1}\frac{1}{(1-q^{6m-5})(1-q^{6m-1})}.\numberthis\label{eq61}\\
\end{align*}
Sometime during the 90s, the second author noted that\\
\begin{align*}
\prod^{\infty}_{m=1}\frac{1}{(1-q^{6m-5})(1-q^{6m-1})}=\prod^{\infty}_{m=1}(1+q^{2m-1}+q^{4m-2}),\numberthis\label{eq62}\\
\end{align*}
where the second product in \eqref{eq62} is the generating function of $A(n)$, the number partitions of $n$ into odd parts repeating no more than twice, and suggested to George Andrews that it would be worthwhile to explore the deeper connections between the equality\\
\begin{align*}
S(n)=A(n).\numberthis\label{eq63}\\
\end{align*}
Andrews dubbed the equality in \eqref{eq63} the \textit{Alladi-Schur Theorem} and
established in \cite{GEA19} the following deep refinement:\\

\textbf{Theorem A:} (Andrews' refinement of the Alladi-Schur theorem)

\textit{Let $A(n;k)$ denote the number of partitions of $n$ into odd parts repeating no more than twice and with exactly $k$ parts}.

\textit{Let $s(n;k)$ denote the number of partitions of $n$ into parts that differ by $\ge 3$, and with no consecutive multiples of $3$, and having exactly $k$ parts, where the even parts are counted twice. Then}\\
\begin{align*}
A(n;k)=s(n;k).\numberthis\label{eq64}\\
\end{align*}

\textbf{Remark 6.1:} What is surprising is that in Theorem A, partitions of the Schur type are classified according to their parity. Thus, the equality \eqref{eq62} has provided a fresh direction for the investigation of Schur partitions. Andrews' proof of Theorem A \cite{GEA19} was $q$-theoretic. In view of the combinatorial elegance of Theorem A, it is natural to ask if there is a combinatorial/bijective proof; such a proof was recently found by Alamoudi in \cite{YA24}, and it turned out to be quite intricate. Some notions in \cite{YA24} share a resemblance with, but are different from, some of the notions in Kur\c{s}ung\"{o}z's important paper \cite{KK21} on Schur's partition theorem.\\

With the combinatorial proof of Theorem A having been found, the following question arises. Can Theorem A be cast in the form of a $q$-hypergeometric key identity? We answer this in the affirmative below.\\

In \cite{KK21}, Kur\c{s}ung\"{o}z obtains, by combinatorial arguments, a series generating function for the Schur partitions, which is different from the series in the Alladi-Gordon key identity; then, by the same combinatorial arguments, he obtains a series generating function for Schur partitions by keeping track of the number of even and odd parts. His result is:\\

\textbf{Theorem K:}
\textit{Let $s(n;m_1,m_0)$ denote the number of Schur partitions of $n$ having $m_1$ odd parts and $m_0$ even parts. Then}\\
\begin{align*}
\sum_{m_1,m_0, n\ge 0}s(n;m_1,m_0)&a^{m_1}b^{m_0}q^n=\sum_{n_{11}, n_{10}, n_{21}, n_{22}\ge 0}\frac{q^{6n^2_{21}-n_{21}+6n^2_{22}+n_{22}+2n^2_{11}-n_{11}+2n^2_{10}}}{(q^2;q^2)_{n_{11}}(q^2;q^2)_{n_{10}}(q^6;q^6)_{n_{21}}(q^6;q^6)_{n_{22}}}\\
&\times \quad q^{12n_{21}n_{22}+6(n_{21}+n_{22})(n_{11}+n_{10})+4n_{11}n_{10}}a^{n_{21}+n_{22}+n_{11}}b^{n_{21}+n_{22}+n_{10}}.\numberthis\label{eq65}\\    
\end{align*}

In deriving Theorem K, Kur\c{s}ung\"{o}z groups the Schur partitions of $n$ into disjoint \textit{pairs}, which are parts that differ by exactly 3 (with a certain convention when there is a maximal chain of $\ell$ parts differing by 3 with $\ell$ is odd), and calls the rest \textit{singletons}. In the above identity, $n_{11}$ (resp. $n_{10}$) is the number of odd (resp. even) singletons, and $n_{21}$ (resp. $n_{22}$) is the number of $1\pmod{3}$ (resp. $2\pmod{3}$) pairs. We now point out that in view of Andrews' refinement of the Alladi-Schur Theorem and Kur\c{s}ung\"{o}z's series representation \eqref{eq65} for the generating function of $s(n;m_1,m_0)$, if we choose\\
\begin{align*}
b=a^2,\numberthis\label{eq66}\\
\end{align*}
then the series in \eqref{eq65} will be equal to the product\\
\begin{align*}
\prod^{\infty}_{m=1}(1+aq^{2m-1}+a^2q^{4m-2}),\numberthis\label{eq67}\\
\end{align*}
and this yields the key identity for Theorem A. That is\\
\begin{align*}
&\sum_{m_1,m_0, n\ge 0}s(n;m_1,m_0)a^{m_1+2m_0}q^n=\sum_{n_{11}, n_{10}, n_{21}, n_{22}\ge 0}\frac{q^{6n^2_{21}-n_{21}+6n^2_{22}+n_{22}+2n^2_{11}-n_{11}+2n^2_{10}}}{(q^2;q^2)_{n_{11}}(q^2;q^2)_{n_{10}}(q^6;q^6)_{n_{21}}(q^6;q^6)_{n_{22}}}\\
&\times \quad q^{12n_{21}n_{22}+6(n_{21}+n_{22})(n_{11}+n_{10})+4n_{11}n_{10}}a^{3n_{21}+3n_{22}+n_{11}+2n_{10}}=\prod^{\infty}_{m=1}(1+aq^{2m-1}+a^2q^{4m-2})\numberthis\label{eq68}\\
\end{align*}
is the \textit{key identity} for Theorem A. This seems to have escaped attention. It is desirable to have a $q$-hypergeometric proof of \eqref{eq68}. Kur\c{s}ung\"{o}z's method has been used by other authors to obtain new series generating functions for various fundamental partition functions. But Theorem K had not been considered in conjunction with Theorem A, which is the reason that the key identity \eqref{eq68} for Theorem A presented here had escaped attention. In \cite{GEA17}, Andrews has expressed the view that the new direction for Schur's theorem presented by Theorem A is deeper and more significant than the classical version of Schur's theorem. This is confirmed by the complexity of the key identity \eqref{eq68}, for which, at the time of this writing, a $q$-hypergeometric proof is not known. There is, however, another series representation for the generating function of $s(n;m_1,m_2)$, due to Andrews-Chern-Li \cite{GEASCZL22}, which is\\
\begin{align*}
\sum_{n,m_1,m_2}s(n;m_1,m_2)&a^{m_1+m_2}b^{m_2}q^n=\sum_{n_1,n_2,n_3\ge 0}\frac{(-1)^{n_3}a^{n_1+n_2+2n_3}b^{n_2+n_3}}{(q^2;q^2)_{n_1}(q^2;q^2)_{n_2}(q^6;q^6)_{n_3}}\\
&\times \quad q^{2n^2_1-n_1+2n^2_2+9n^2_3+2n_1n_2+6n_1n_3+6n_2n_3}.\numberthis\label{eq69}\\
\end{align*}
When one sets $a=b$ in \eqref{eq69}, which means that the even parts are counted twice, then one can set the resulting expression equal to the product on the right in \eqref{eq68}. Thus, we have the identity\\

\begin{align*}
\sum_{n_1,n_2,n_3\ge 0}\frac{(-1)^{n_3}a^{n_1+2n_2+3n_3}}{(q^2;q^2)_{n_1}(q^2;q^2)_{n_2}(q^6;q^6)_{n_3}}&\times \quad q^{2n^2_1-n_1+2n^2_2+9n^2_2+2n_1n_2+6n_1n_3+6n_2n_3}\\
&=\prod^{\infty}_{m=1}(1+aq^{2m-1}+a^2q^{4m-2}).\numberthis\label{eq610}\\
\end{align*}
Andrews-Chern-Li provide two proofs of \eqref{eq610}, one $q$-hypergeometric, and another which is computer-aided. But it is to be noted that in \eqref{eq69} and \eqref{eq610}, there is the factor $(-1)^{n_3}$, and so it is not transparent that the coefficients in the power series expansion are non-negative. On the other hand, it is transparent that the coefficients of the series on the left in \eqref{eq68} are all non-negative.\\

\section{The Capparelli theorems and the key-identity}\label{s7}
Through a study of vertex operators in Lie algebras, Capparelli \cite{SC93} conjectured the following partition theorem:\\

\textbf{Theorem C:}
\textit{Let $C^*(n)$ denote the number of partitions of $n$ into parts $\equiv\pm 2, \pm 3\pmod{12}$.}

\textit{Let $D(n)$ denote the number of partitions of $n$ into parts $>1$ with minimal difference $\ge 2$, where the difference is $\ge 4$ unless consecutive parts are multiples of $3$ or add up to a multiple of $6$. Then}\\
\begin{align*}
C^*(n)=D(n).\\
\end{align*}

The first proof of Theorem C was due to Andrews \cite{GEA94} by the use of generating functions. Subsequently, Alladi-Andrews-Gordon \cite{KAGEABG95} noticed that if $C^*(n)$ is replaced by the equivalent partition function $C(n)$, which is the number of partitions of $n$ into distinct parts $\equiv 2, 3, 4\, \text{or}\, 6 \pmod{6}$, then there is a three-parameter refinement, namely:\\

\textbf{Theorem C-R:}
\textit{Let $C(n;i,j,k)$ denote the number of partitions of the type enumerated by $C(n)$, with the added restriction that there are precisely $i$ parts $\equiv 4\pmod{6}$, $j$ parts $\equiv 2\pmod{6}$, and of those $\equiv 0\pmod{3}$, exactly $k$ are $>3(i+j)$.}

\textit{Let $D(n;i,j,k)$ denote the number of partitions of the type enumerated by $D(n)$ with the additional restriction that there are precisely $i$ parts $\equiv 1\pmod{3}$, $j$ parts $\equiv 2\pmod{3}$, and $k$ parts $\equiv 0\pmod{3}$. Then}\\
\begin{align*}
C(n;i,j,k)=D(n;i,j,k).\\
\end{align*}

Alladi-Andrews-Gordon \cite{KAGEABG95} established a generalization of Theorem C-R by the \textit{method of weighted words} (which was initiated in Alladi-Gordon \cite{KABG93MM} to establish a generalization of Schur's theorem), and viewed this generalized theorem as the combinatorial interpretation of the following \textit{key identity}:\\
\begin{align*}
&\sum_{i,j}\frac{a^ib^jq^{2T_i+2T_j}(-q)_{i+j}(-cq^{i+j+1})_{\infty}}{(q^2;q^2)_i(q^2;q^2)_j}\\
&=\sum_{i,j,k}\frac{a^ib^jc^kq^{2T_i+2T_j+T_k+(i+j)k}}{(q)_{i+j+k}}{{i+j+k} \brack k}_q{{i+j} \brack i}_{q^2}.\numberthis\label{eq71}\\
\end{align*}

The sum on the right in \eqref{eq71} could be rewritten as\\
\begin{align*}
&\sum_{i,j,k}\frac{a^ib^jc^kq^{2T_i+2T_j+T_k+(i+j)k}}{(q)_{i+j+k}}\frac{(q)_{i+j+k}}{(q)_{i+j}(q)_k}\frac{(q^2;q^2)_{i+j}}{(q^2;q^2)_i(q^2;q^2)_j}\\
&=\sum_{i,j}\frac{a^ib^jq^{2T_i+2T_j}(-q)_{i+j}}{(q^2;q^2)_i(q^2;q^2)_j}\sum_k\frac{c^kq^{T_k+(i+j)k}}{(q)_k}\\
&=\sum_{i,j}\frac{a^ib^jq^{2T_i+2T_j}(-q)_{i+j}(-cq^{i+j+1})_{\infty}}{(q^2;q^2)_i(q^2;q^2)_j}.\numberthis\label{eq72}\\
\end{align*}
If we take $c=1$, then the term in \eqref{eq72} becomes\\
\begin{align*}
(-q)_{\infty}\sum_{i,j}\frac{a^ib^jq^{2T_i+2T_j}}{(q^2;q^2)_i(q^2;q^2)_j}=(-aq^2;q^2)_{\infty}(-bq^2;q^2)_{\infty}(-q)_{\infty}.\numberthis\label{eq73}\\    
\end{align*}
If we make the replacements\\
\begin{align*}
q\mapsto q^3, \quad a\mapsto q^{-2}, \quad b\mapsto q^{-4}\numberthis\label{eq74}\\
\end{align*}
in \eqref{eq74}, we get the generating function of $C(n)$ in Theorem C-R.\\

Capparelli \cite{SC93} had stated another partition conjecture in the form $A^*(n)=B(n)$. The difference between the conditions defining $B(n)$ and $D(n)$ is that among the partitions enumerated by $B(n)$, the integer $2$ should not occur as a part, but $1$ is allowed as part. The generating function of $A^*(n)$ is a product that is more complicated than the product generating function of $C^*(n)$. However, it turns out that $A^*(n)=A(n)$, where $A(n)$ is the number of partitions of $n$ into distinct parts $\equiv 1,3,5, \textit{or } 6\pmod{6}$. So, this second conjecture is equivalent to $A(n)=B(n)$. This conjecture can be proved by applying the transformations\\
\begin{align*}
q\mapsto q^3, \quad a\mapsto q^{-5}, \quad b\mapsto q^{-1},\numberthis\label{eq75}\\
\end{align*}
and by combinatorially interpreting the resulting $q$-hypergeometric identity. So what we stress here is that in \cite{KAGEABG95}, by considering partitions into \textit{distinct} parts in certain residue classes modulo 6, instead of partitions into parts in certain residue classes modulo 12 (parts that could repeat), not only is the second partition theorem of Capparelli cast in a more elegant form, but also that such a reformulation is capable of refinement and generalization. The idea that reformulating partitions into certain distinct parts is capable of refinements was initiated in Alladi-Gordon's treatment of Schur's theorem \cite{KABG93MM} and indeed that was instrumental in developing the \textit{method of weighted words} which is widely applicable.\\

It is to be noted that when we set $a=0$, the generalized Capparelli product on the right-hand side of \eqref{eq73} reduces to\\
\begin{align*}
(-bq;q^2)_{\infty}(-q)_{\infty},\numberthis\label{eq76}\\    
\end{align*}
which is the product for Lebesgue's identity. In view of this link between the generalized Capparelli identity and Lebesgue's identity, it is natural to ask whether a finite version of the Capparelli key identity can be obtained from the Transformation Formula (Lemma 3.3) in Section \ref{s3} just as we obtained the finite version of Lebesgue's identity from Lemma 3.3. The answer is YES. Indeed, Alladi (1994 unpublished) obtained the following from Lemma 3.3: \textit{For positive integers $m,n$, we have}\\
\begin{align*}
&\sum^m_{j=0}\sum^n_{i=0}(bc)^ia^jq^{2T_i+2T_j}{m \brack i}_q{n \brack j}_{q^2}(-cq^{i+1})_{\ell+j}\\
&=\sum^{m+n+\ell}_{N=0}c^Nq^{T_N}\left(\sum^n_{j=0}\sum^m_{i=0}b^ia^jq^{T_i+2T_j}{m \brack i}_q{n \brack j}_{q^2}{{\ell+j} \brack {N-i}}_q\right).\numberthis\label{eq77}\\
\end{align*}
If we let $m,n,\ell\to\infty$ in \eqref{eq77}, then we get\\
\begin{align*}
&\sum^{\infty}_{i=0}\sum^{\infty}_{j=0}\frac{(bc)^iia^jq^{2T_i+2T_j}(-q)_i(-cq^{i+1})_{\infty}}{(q^2;q^2)_i(q^2;q^2)_j}\\
&=(-aq^2;q^2)_{\infty}\sum^{\infty}_{N=0}\frac{c^Nq^{T_N}(-bq)_N}{(q)_N}.\numberthis\label{eq78}\\
\end{align*}
Identities \eqref{eq77} and \eqref{eq78} are different from the Capparelli key identity \eqref{eq71} in the sense that in \eqref{eq71} we have $(-q)_{i+j}(-cq^{i+j+1})_{\infty}$, whereas in \eqref{eq78} we have $(-q)_i(-cq^{i+1})_{\infty}$. However, when $c=1$ both versions are the same, and the expression on \eqref{eq78} becomes the product\\
\begin{align*}
(-aq^2;q^2)_{\infty}(-bq^2;q^2)_{\infty}(-q)_{\infty}.\numberthis\label{eq79}\\
\end{align*}
But there is a general version of \eqref{eq71}, which specializes into a polynomial identity, with a free parameter $c$ where the decomposition involves $(-q)_{i+j}$ as noted by Alamoudi: \textit{For integers $M_1, M_2, L$, we have}\\
\begin{align*}
&\sum_{j,i\geq 0}a^ib^jq^{2T_i+2T_j}{{M_1} \brack i}_{q^2}{{M_2}\brack j}_{q^2}(-q)_{i+j}(-cq^{i+j+1})_{L-i-j}\\
&=\sum_{i,j,k\geq 0}a^ib^jc^kq^{2T_i+2T_j+T_k+(i+j)k} \frac{(q^{L-i-j-k+1})_k(q^{2M_1-2i+2};q^2)_i(q^{M_2-2j+2};q^2)_j}{(q)_{i+j+k}}{{i+j+k} \brack k}{{i+j} \brack i}_{q^2}\\
\numberthis\label{eq710}\\
\end{align*}
\textit{In particular, for non-negative integers $M_1, M_2, L$, with $L\geq M_1+M_2$ we have}\\
\begin{align*}
&\sum^{M_1}_{i=0}\sum^{M_2}_{j=0}a^ib^jq^{2T_i+2T_j}{{M_1} \brack i}_{q^2}{{M_2}\brack j}_{q^2}(-q)_{i+j}(-cq^{i+j+1})_{L-i-j}\\
&=\sum_{\substack{i,j,k\geq 0\\i\le M_1,\, j\le M_2,\,i+j+k\leq L}}a^ib^jc^kq^{2T_i+2T_j+T_k+(i+j)k}{L \brack {i+j+k}}{{i+j+k} \brack k}{{i+j} \brack i}_{q^2}\\
&\qquad\qquad\qquad\qquad\qquad\times\quad\frac{(q^{2M_1-2i+2};q^2)_i(q^{M_2-2j+2};q^2)_j}{(q^{L-i-j+1})_{i+j}}.\numberthis\label{eq711}\\
\end{align*}
After performing various cancellations, the right-hand sides of \eqref{eq710} and \eqref{eq711} become \\
\begin{align*}
\sum_{i,j,k\geq 0}a^ib^jc^kq^{2T_i+2T_j+T_k+(i+j)k}{{M_1} \brack i}_{q^2}{{M_2}\brack j}_{q^2}{{L-i-j} \brack k}(-q)_{i+j}.\numberthis\label{eq712}\\
\end{align*}
If we sum the inner sum over $k$ and use\\
\begin{align*}
\sum_{k\geq0}c^kq^{T_k+(i+j)k}{{L-i-j}\brack k}=(-cq^{i+j+1})_{L-i-j},\numberthis\label{eq713}\\
\end{align*}
we get the left-hand sides of \eqref{eq710} and \eqref{eq711}. If we let $M_1, M_2, L\to\infty$, we get the (infinite version) Capparelli key identity given by \eqref{eq71} and \eqref{eq72}.\\

\textbf{Remark 7.1:} We have given the intermediate identities \eqref{eq710} and \eqref{eq711} because they maintain the form of \eqref{eq71}, which is the key identity of the original Alladi-Andrews-Gordon three-parameter refinement of the Capparelli partition theorem. The combinatorial significance of the form in \eqref{eq71} is that it highlights the generating function \cite[Eq. (5.5)]{KAGEABG95} which counts the relevant minimal partitions.\\

\textbf{Remark 7.2:} Given the generality of \eqref{eq710}, we can apply the same bound shifting technique used to obtain \eqref{eq56} from Theorem 5.1 and get other finite analogs of \eqref{eq71}. Specifically, by setting the coefficients of $a^ib^j$ equal, we get
\begin{align*}
&a^ib^jq^{2T_i+2T_j}{{M_1} \brack i}_{q^2}{{M_2}\brack j}_{q^2}(-q)_{i+j}(-cq^{i+j+1})_{L-i-j}\\
&=\sum_{k\geq 0}a^ib^jc^kq^{2T_i+2T_j+T_k+(i+j)k} \frac{(q^{L-i-j-k+1})_k(q^{2M_1-2i+2};q^2)_i(q^{M_2-2j+2};q^2)_j}{(q)_{i+j+k}}{{i+j+k} \brack k}{{i+j} \brack i}_{q^2}.\\
\numberthis\label{eq714}
\end{align*}
The above equation is valid for any choice of integers $M_1,M_2,L,i,j$ with $i,j\geq 0$. Setting $M_1=L-j$, $M_2=L$, canceling using $(q^{a-j-i+1})_{i+j}=(q^{a-j-i+1})_{i}(q^{a-j+1})_{j}$ for any $a\in\mathbb{Z}$, and summing for all $i,j\geq 0$ gives, for any $L\in\mathbb{Z}$,
\begin{align*}
&\sum_{i,j\geq0}a^ib^jq^{2T_i+2T_j}{{L-j} \brack i}_{q^2}{{L}\brack j}_{q^2}(-q)_{i+j}(-cq^{i+j+1})_{L-i-j}\\
&=\sum_{i,j,k\geq 0}a^ib^jc^kq^{2T_i+2T_j+T_k+(i+j)k}{L \brack {i+j+k}}{{i+j+k} \brack k}{{i+j} \brack i}_{q^2}(-q^{L-i-j+1})_{i+j}.\\
\numberthis\label{eq715}
\end{align*}\\
The point is in \eqref{eq714}, the integers $M_1,M_2,i,j,$ and $L$ are fixed and only $k$ is being summed over. Thus, many substitutions can be made, perhaps sending $M_1$ to something that depends on $j$, and the expression remains valid, after which one can sum over $i,j\geq0$. This way, one may obtain a myriad of other finite analogs of \eqref{eq71}. However, \textbf{the same liberty does not necessarily extend to shifting by expressions involving $k$} since $k$ is still bound by the summation in \eqref{eq714} and is not free as $i\geq 0$ and $j\geq 0$ are in \eqref{eq714}.\\

In closing this section, we note that when we set $a=0$ in \eqref{eq79}, the Capparelli product generating function reduces to the product generating function for Lebesgue's identity, and if we further set $b=0$, we get the product generating function for Euler's theorem on partitions into distinct parts. This leads us to an infinite hierarchy of identities observed by Alladi (unpublished) in 1994, of which the first three cases are those of Euler, Lebesgue, and Capparelli, in that order. We present this in the next section, along with a new polynomial version of this infinite hierarchy due to Alamoudi.\\

\section{An infinite hierarchy of $q$-hypergeometric identities}\label{s8}
In 1994, Alladi observed that for each non-negative integer $r$, there is the identity\\
\begin{align*}
&\sum_{\nu_1, \nu_2, \cdots, \nu_r, k\ge 0}\frac{a^{\nu_1}_1a^{\nu_2}_2\cdots a^{\nu_r}_rc^kq^{2T_{\nu_1}+2T_{\nu_2}\cdots +2T_{\nu_r}+T_k+k(\nu_1+\nu_2+\cdots \nu_k)}(-q)_{\nu_1+\nu_2+\cdots\nu_k}}{(q^2;q^2)_{\nu_1}(q^2;q^2)_{\nu_2}\cdots (q^2;q^2)_{\nu_r}(q)_k}\\
&=\sum_{\nu_1, \nu_2, \cdots, \nu_r\ge 0}\frac{a^{\nu_1}_1a^{\nu_2}_2\cdots a^{\nu_r}_rq^{2T_{\nu_1}+2T_{\nu_2}\cdots +2T_{\nu_r}}(-q)_{\nu_1+\nu_2+\cdots\nu_r}(-cq^{\nu_1+\nu_2+\cdots +\nu_r+1})_{\infty}}{(q^2;q^2)_{\nu_1}(q^2;q^2)_{\nu_2}\cdots (q^2;q^2)_{\nu_r}}.\numberthis\label{eq81}\\
\end{align*}
Identity \eqref{eq81} is easily proved by summing over $k$ the inner sum on the left-hand side. When $c=1$, identity \eqref{eq81} becomes\\
\begin{align*}
\sum_{\nu_1, \nu_2, \cdots, \nu_r, k\ge 0}&\frac{a^{\nu_1}_1a^{\nu_2}_2\cdots a^{\nu_r}_rq^{2T_{\nu_1}+2T_{\nu_2}\cdots +2T_{\nu_r}+T_k+k(\nu_1+\nu_2+\cdots \nu_k)}(-q)_{\nu_1+\nu_2+\cdots\nu_k}}{(q^2;q^2)_{\nu_1}(q^2;q^2)_{\nu_2}\cdots (q^2;q^2)_{\nu_r}(q)_k}\\
&=(-a_1q^2;q^2)_{\infty}(-a_2q^2;q^2)_{\infty}\cdots(a_rq^2;q^2)_{\infty}(-q)_\infty.\numberthis\label{eq82}\\    
\end{align*}
Note that the case $r=0$ in \eqref{eq82} is Euler's identity for partitions into distinct parts, the case $r=1$ gives the product for Lebesgue's identity, and the case $r=2$ yields the product for the generalized Capparelli identity \eqref{eq73}. This was the motivation to come up with this infinite hierarchy:\\
\begin{align*}
\text{Euler}\, (r=0) \quad \to \quad \text{Lebesgue}\, (r=1)\quad \to \text{Capparelli}\, (r=2)\quad \to\quad\cdots .\numberthis\label{eq83}\\
\end{align*}

To get regular partition theorems (not weighted ones) from \eqref{eq81} and \eqref{eq82}, the minimal dilation is\\
\begin{align*}
q\mapsto q^{r+1},\numberthis\label{eq84}\\
\end{align*}
and this is the optimal \textit{dilation}. There are several choices of residue classes modulo $r+1=m$ for this dilation: Let $j_1, j_3, \cdots, j_r$ be incongruent modulo $m$, with $0<j_i<2m$, $j_i\ne m$. Consider the \textit{translations}\\
\begin{align*}
a_i\mapsto a_iq^{j_i-m}\quad\text{for}\quad i=1,2,\cdots r.\numberthis\label{eq85}\\
\end{align*}
We will now turn our attention to the corresponding partition theorem. However, we must first discuss the important notion of level parities. We begin by alerting the reader to the following convention we adopt in the sequel.\\

In the sequel, we adopt the following convention. Whenever we say that a part is $\equiv j\pmod{M}$, we mean a part of the form $j+\lambda M$, with the integer $\lambda\ge 0$. In particular, the part is $\ge j$. For example, if a part is $\equiv 4\pmod{3}$, that part must be one of $4,7,10, \cdots$.\\

For each part $p\equiv j\pmod{M}$, such that $M\nmid j$, we say that $p$ has \textbf{odd level parity} as a\\ $ j\pmod{M}$ part if $p=j+\lambda M$ with $\lambda$ odd. Otherwise, we say $p$ has \textbf{even level parity} as a $ j\pmod{M}$ part.\\ 

\textbf{Remark 8.1:}
In the sequel, we will simply say that a part $p$ has odd/even level parity and omit the prepositional phrase ``as a $ j\pmod{M}$ part" as it will be clear from context. Nonetheless, we would like to alert the reader to the following subtlety. In general, the level parity of a part $p$ depends on the specific choice of $j$ and $M$. For example, $7$ has odd level parity if $j=2$ and $M=5$, but even level parity if $j=1$ and $M=3$. Moreover, in view of the above convention, even if $j\equiv j'\pmod{M}$ in the traditional sense, $j$ and $j'$ can define different level parities for $p$. For example, $7$ has odd level parity if $j=4$ and $M=3$ but even level parity if $j=1$ and $M=3$, even though $4\equiv1\pmod{3}$.\\
\newpage

\textbf{Theorem H:}
\textit{Let $A(n;\nu_1, \nu_2, \cdots, \nu_r; 2m)$ denote the number of partitions of $n$ into $\nu_i$ distinct parts $\equiv j_i\pmod{2m}$ and distinct parts $\equiv0\pmod{m}$, for $i=1,2,\cdots, r$.}

\textit{Let $B(n; \nu_1,\nu_2,\cdots,\nu_r; m)$ denote the number of partitions of $n$ into $\nu_i$ distinct parts\\ $\equiv j_i\pmod{m}$ for $i=1,2,\cdots, r$ such that the difference between two parts of different level parities is $>m$, and distinct parts $\equiv 0\pmod{m}$ each $>(j_1+j_2+\cdots+j_r)m$.}

\textit{Then, we have}
\begin{align*}
A(n;\nu_1, \nu_2, \cdots, \nu_r; 2m)=B(n; \nu_1,\nu_2,\cdots,\nu_r;m).\numberthis\label{eq86}\\    
\end{align*}

\textbf{Remark 8.2}: In the combinatorial proofs of Schur's theorem in Alladi-Gordon \cite{KABG93MM}, and of Capparelli's theorem in Alladi-Andrews-Gordon \cite{KAGEABG95}, a combinatorial method due to Bressoud is followed. This method involves two stages: Stage 1 is an \textit{embedding} (of the Ferrers conjugate of the distinct parts $\equiv 0\pmod{m}$ which are $\le (\nu_1+\nu_2+\cdots+\nu_r)m$ into the Ferrers graph of the partitions parts into distinct parts $\equiv j_i\pmod{2m}$ for $i=1,2,\cdots,r$. This is followed by Stage 2, which is a rearrangement, which we call the \textit{Bressoud rearrangement}, and this is more complicated. It is only after the rearrangement that we get a partition satisfying certain difference conditions. The partition function $A$ in (8.6) is at the product level, and this is the start of the combinatorial construction. The partition function $B$ is at the embedding stage of the combinatorial proof. What we require is the partition function $C$ (whose parts satisfy difference conditions based on the residue classes), which results AFTER the Bressoud rearrangement is completed. We expect to present the function $C$ after working out the details in a subsequent paper \cite{YAKA25} devoted to the combinatorics of various identities that have been presented here. However, we will present C in a certain special case (see Theorem C5). \\

\textbf{Remark 8.3}: We note the following with regards to Theorem H.
\begin{enumerate}[(i)]
    \item It is to be noted that for the partition function $A$ in Theorem H, the modulus is $2m$, whereas for the partition function $B$, the modulus is $m$. This halving of the modulus for $B$ is due to the embedding.\\
    \item  Due to observed combinatorial insights and intricacies, the authors have refined the statement of Theorem H, as well as Theorem C5 presented later in this section, in this version of the manuscript. This includes the introduction of the notion of level parities. Regarding this notion, it is important to note that although the difference between two parts of different level parities is $>m$, in view of the definition of level parity, consecutive multiples of $m$ are allowed for Theorem H. 
\end{enumerate}

\medskip 

In the case $r=1$, the dilation gives the Little G\"ollnitz theorems, not the Lebesgue identity which corresponds to the undilated case.\\

For $r=1$, the function $C$ is as in Corollary $2^A$ of Alladi-Gordon \cite{KABG93JCTA}. There, the function $B$ is not mentioned. It would be $B(n;j,2)=$ \textit{the number of partitions into $j$ distinct and non-consecutive odd parts, and even parts, each $>2j$}.\\

In the case $r=2$, which corresponds to Capparelli's theorem, $C$ and $A$ are known, but $B$ is not mentioned. It will be $B(n;\nu_1,\nu_2;3)=$ \textit{the number of partitions of $n$ into distinct $\nu_1$ parts $\equiv 4\pmod{3}$ and distinct $\nu_2$ parts $\equiv 2\pmod{3}$ such that the gap between parts of different parities is $>3$, and distinct parts $\equiv 0\pmod{3}$ each $>3(\nu_1+\nu_2)$}.\\ 

Actually, there are two theorems of Capparelli (which correspond to the case $r=2$). For the first theorem, the function $A$ has generating function\\
\begin{align*}
(-a_1q^2;q^6)_{\infty}(-a_2q^4;q^6)_{\infty}(-q^3;q^3)_{\infty}.\numberthis\label{eq87}\\
\end{align*}
The second theorem of Capparelli corresponds to the generating function of $A(n)$ being\\
\begin{align*}
(-a_1q;q^6)_{\infty}(-a_2q^5;q^6)_{\infty}(-q^3;q^3)_{\infty}.\numberthis\label{eq88}\\
\end{align*}
When we set $a_1=a_2=1$, the product in \eqref{eq87} becomes\\
\begin{align*}
\frac{1}{(q^3;q^6)_{\infty}(q^2;q^{12})_{\infty}(q^{10};q^{12})_{\infty}}=\frac{1}{(q^3;q^{12})_{\infty}(q^9;q^{12})_{\infty}(q^2;q^{12})_{\infty}(q^{10};q^{12})_{\infty}}.\numberthis\label{eq89}\\    
\end{align*}
It is in the unrefined form as in \eqref{eq89} that Capparelli stated his conjecture. It was in Alladi-Andrews-Gordon \cite{KAGEABG95} that the product in \eqref{eq89} was replaced by \eqref{eq87}, which had the advantage that the theorem could be refined by introducing parameters $a_1, a_2$. When we set $a_1=a_2=1$ in \eqref{eq88}, it does not yield a product as nice as in \eqref{eq89}, but more complicated, and it is such a product that Capparelli stated. Again, the advantage of considering distinct parts, as pointed out in \cite{KAGEABG95}, is that Capparelli's second theorem can be more neatly expressed in terms of the product \eqref{eq88}.  So the question that arises is whether, for $r\ge 3$,  there are such special choices of the residue classes modulo $2m=2(r+1)$, where the product in terms of distinct parts determined modulo $m$ transforms neatly into a product where the modulus is $2m$, and the parts can repeat as in \eqref{eq89}. Such a phenomenon does not occur when $r=3$, but does occur for all even $r\ge 4$. We illustrate this with an example of $r=4$ (so $m=r+1=5$). Thus, as per Theorem L, we need to choose four residue classes $\mod 10$, which are incongruent $\mod 5$ and unequal to 0 or 5 $\pmod{10}$. The ideal choice is $2,4,6,8\pmod{10}$. This then gives\\
\begin{align*}
&(-q^2;q^{10})_{\infty}(-q^4;q^{10})_{\infty}(-q^6;q^{10})_{\infty}(-q^8;q^{10})_{\infty}(-q^5;q^5)_{\infty}\\
&=\frac{1}{(q^5;q^{10})_{\infty}(q^2;q^{20})_{\infty}(q^6;q^{20})_{\infty}(q^{14};q^{20})_{\infty}(q^{18};q^{20})_{\infty}},\numberthis\label{eq810}\\
\end{align*}
which is an ideal extension of Capparelli's first theorem to the level $r=4$. Thus, the combinatorics underlying this infinite hierarchy is fascinating, and we hope to discuss this in a subsequent paper \cite{YAKA25}. But for this paper, we will only state the Capparelli theorem in the higher case $r=4$ ($m=r+1=5$), which we can prove combinatorially by the method given in Alladi-Andrews-Gordon \cite{KAGEABG95}. \\


\textbf{Theorem C5:}
\textit{Let $A(n)$ denote the number of partitions of $n$ into distinct parts $\equiv 0,2,4,5,\\6,\, \text{or} \, 8\pmod{10}$.}

\textit{Let $C(n)$ denote the number of partitions of $n$ into parts $\equiv 2, 5, 6, 14, 15, \, \text{or} \, 18 \pmod{20}$.}

\textit{Let $D(n)$ denote the number of partitions $\lambda=\lambda_1+\dots \lambda_\ell$, written in non-increasing order, with distinct parts not equal to $1$ or $3$, such that:}\begin{enumerate}[(i)]
    \item \textit{The gap is $\geq 5$ if two consecutive parts have different parities or one of them is a multiple of $5$, with gap $5$, only allowed for consecutive multiples of $5$.}\\
    \item \textit{Whenever $5|\lambda_i$, we have $\lambda_i-\lambda_{i+j}\geq 5j$ for $1\leq i\leq i+j\leq \ell$.}
\end{enumerate}
\textit{Then,}
\begin{align*}
A(n)=C(n)=D(n).\\
\end{align*}


\textbf{Remark 8.4:} In \cite{YAKA25}, we will consider a general form of the second condition for $D(n)$. Namely, for modulus $m$ we consider partitions with the property that whenever $m|\lambda_i$, we have $\lambda_i-\lambda_{i+j}\geq mj$ for $1\leq i\leq i+j\leq \ell$. We will also consider a dual variant of this notion. These considerations will be in the context of higher analogs of $D(n)$ featured in Theorem C5.  Conversely, for lower $m$, namely for Capparalli (i.e. $m=3$), we note that condition (ii) is redundant.\\

The combinatorial proof actually gives a four-parameter refinement of the equality $A(n)=D(n)$ in which we can keep track of the parts $\equiv 1,2,3,4\pmod{5}$. Note that for the partition function $A(n)$ in Theorem C5, the parts $\equiv 1,2,3,4\pmod{5}$ are actually in the form $6,2,8,4\pmod{10}$. All this will be presented in a subsequent paper along with the combinatorics of the infinite hierarchy.\\

For this manuscript, however, we will conclude with a discussion of the polynomial infinite hierarchy. In particular, Alamoudi has recently noticed that the following polynomial (finite) identities correspond to \eqref{eq81} and \eqref{eq82}. For simplicity of the expressions, put $N=\nu_1+\nu_2+\cdots+\nu_r$. Then, for integers $M_1, M_2, \dots,M_r, L$, we have\\
\begin{align*}
 \sum_{\nu_1,\dots,\nu_r,k\geq 0}&a^{\nu_1}_1 q^{T_{\nu_1}}_1a^{\nu_2}_2q^{T_{\nu_2}}_2\cdots a^{\nu_r}_r q^{T_{\nu_r}}_r q^{T_k + kN}{{M_1}\brack {\nu_1}}_{q_{1}} {{M_2}\brack {\nu_2}}_{q_{2}}\cdots {{M_r}\brack {\nu_r}}_{q_{r}}{{L-N}\brack k} (-q)_{N}\\
&=(-a_1q_{1};q_1)_{M_1}(-a_2q_2;q_2)_{M_2}\cdots(-a_rq_r;q_r)_{M_r}(-q)_{L}.\numberthis\label{eq811}\\
\end{align*}
Notice that when $M_1, M_2, \dots, M_r, L$ are non-negative, with $L\geq M_1+M_2+ \dots+M_r$, the sum on the left is equivalent to the conditions that\\
\begin{align*}
0\le \nu_i\le M_i \quad \text{for} \quad i=1,2,\cdots, r, \quad \text{and}
\quad N+k\le L.\numberthis\label{eq812}\\
\end{align*}

To prove \eqref{eq811}, sum the inner sum over $k$ and use\\
\begin{align*}
\sum_{k\geq0}c^kq^{T_k+Nk}{{L-N}\brack k}=(-cq^{N+1})_{L-N},\numberthis\label{eq813}
\end{align*}
and
\begin{align*}
(-a_iq_i;q_i)_{M_i}=\sum_{\nu_i\geq0}a_iq_{i}^{T_{\nu_i}}{{M_i}\brack \nu}_{q_i} \quad \text{for} \quad i=1,2,\cdots, r.\numberthis\label{eq814}\\
\end{align*}
When $q_1=\cdots=q_r=q^2$, we obtain\\
\begin{align*}
\sum_{\nu_1,\dots,\nu_r\geq 0}&a^{\nu_1}_1a^{\nu_2}_2\cdots a^{\nu_r}_rc^kq^{2(T_{\nu_1}+T_{\nu_2}+\cdots+T_{\nu_r})+T_k+Nk}{{M_1}\brack {\nu_1}}_{q^2}{{M_2}\brack {\nu_2}}_{q^2}\cdots {{M_r}\brack {\nu_r}}_{q^2}{{L-N}\brack k}(-q)_{N}\\
&=(-a_1q^2;q^2)_{M_1}(-a_2q^2;q^2)_{M_2}\cdots(-a_rq^2;q^2)_{M_r}(-q)_{L}.\numberthis\label{eq815}\\
\end{align*}

Moreover, when $M_i$ for $i=1,2,\cdots,r$ and $L$ all $\to\infty$, \eqref{eq811} reduces to \eqref{eq82}. On the other hand, for a polynomial analog of \eqref{eq81} maintaining the form of \eqref{eq71}, we first write the LHS of \eqref{eq81} in the form of \eqref{eq71} to get\\
\begin{align*}\label{eq81b}
\sum_{\nu_1, \nu_2, \cdots, \nu_r, k\ge 0}&\frac{a^{\nu_1}_1a^{\nu_2}_2\cdots a^{\nu_r}_rc^kq^{2T_{\nu_1}+2T_{\nu_2}\cdots +2T_{\nu_r}+T_k+k(\nu_1+\nu_2+\cdots \nu_r)}}{(q)_{\nu_1+\nu_2+\cdots\nu_r+k}}{{\nu_1+\cdots+\nu_r}\brack {\nu_1, \cdots, \nu_r}}_{q^2}\\
&\times \quad {{\nu_1+\cdots+\nu_r+k}\brack {\nu_1+\cdots+\nu_r,k}}_{q}\\
&=\sum_{\nu_1, \nu_2, \cdots, \nu_r\ge 0}\frac{a^{\nu_1}_1a^{\nu_2}_2\cdots a^{\nu_r}_rq^{2T_{\nu_1}+2T_{\nu_2}\cdots +2T_{\nu_r}}(-q)_{\nu_1+\nu_2+\cdots\nu_r}(-cq^{\nu_1+\nu_2+\cdots +\nu_r+1})_{\infty}}{(q^2;q^2)_{\nu_1}(q^2;q^2)_{\nu_2}\cdots (q^2;q^2)_{\nu_r}}.\tag{8.1b}\\
\end{align*}

As before, for simplicity of the expressions, put $N=\nu_1+\nu_2+\cdots+\nu_r$. Then\\
\begin{align*}
&\sum_{\nu_1, \nu_2, \cdots, \nu_r, k\ge 0}a^{\nu_1}_1a^{\nu_2}_2\cdots a^{\nu_r}_rc^kq^{2(T_{\nu_1}+T_{\nu_2}+\cdots+T_{\nu_r})+T_k+k(\nu_1+\nu_2+\cdots \nu_r)}\frac{(q^{L-N-k+1})_k}{(q)_{N+k}}{{N}\brack {\nu_1, \cdots, \nu_r}}_{q^2}\\
&\qquad\qquad\qquad\qquad\times \quad 
{{N+k}\brack {N,k}}_{q}\prod_{i=1}^r(q^{2M_i-2\nu_i+2};q^2)_{\nu_i}\\
&=\sum_{\nu_1, \nu_2, \cdots, \nu_r\ge 0}a^{\nu_1}_1a^{\nu_2}_2\cdots a^{\nu_r}_rq^{2(T_{\nu_1}+T_{\nu_2}+\cdots+T_{\nu_r})}{{M_1}\brack {\nu_1}}_{q^2}{{M_2}\brack {\nu_2}}_{q^2}\cdots {{M_r}\brack {\nu_r}}_{q^2}(-q)_N(-cq^{N+1})_{L-N}.\numberthis\label{eq816}\\
\end{align*}
In particular, for non-negative integers $M_1, \dots,M_r, L$, with $L\geq M_1+ \dots+M_r$ we have
\begin{align*}
&\sum_{\substack{\nu_1, \nu_2, \cdots, \nu_r, k\geq 0\\ \nu_1 \leq M_1, \cdots, \nu_r \leq M_r, N+k\leq L}}a^{\nu_1}_1a^{\nu_2}_2\cdots a^{\nu_r}_rc^kq^{2(\sum_{i=1}^{r}T_{\nu_i})+T_k+k(\nu_1+\nu_2+\cdots \nu_r)}{{L}\brack {N+k}}_{q}{{N}\brack {\nu_1, \cdots, \nu_r}}_{q^2}\\
&\qquad\qquad\qquad\qquad\times \quad {{N+k}\brack {N,k}}_{q}\frac{\prod_{i=1}^r(q^{2M_i-2\nu_i+2};q^2)_{\nu_i}}{(q^{L-N+1})_{N}}\\
&=\sum_{\nu_1= 0}^{M_1}\sum_{\nu_2= 0}^{M_2}\dots\sum_{\nu_r= 0}^{M_r}a^{\nu_1}_1a^{\nu_2}_2\cdots a^{\nu_r}_rq^{2(\sum_{i=1}^{r}T_{\nu_i})}{{M_1}\brack {\nu_1}}_{q^2}{{M_2}\brack {\nu_2}}_{q^2}\cdots {{M_r}\brack {\nu_r}}_{q^2}(-q)_N(-cq^{N+1})_{L-N}.\numberthis\label{eq817}
\end{align*}\medskip

In the same way as before, when $M_i$ for $i=1,2,\cdots,r$, and $L$ all $\to\infty$, \eqref{eq816} and \eqref{eq817} give \eqref{eq81b}. We note that just like \eqref{eq710}, the LHS of \eqref{eq816} could be written in the form of the LHS of \eqref{eq811} or \eqref{eq815}, then the identity is evident. However, we emphasize this form because it maintains its resemblance to the original \eqref{eq71}. Furthermore, the form\footnote{See Remark 8.7.} of \eqref{eq811} does not exploit the pivotal role the substitution $q_1=\cdots=q_r=q^2$ plays in giving this form. For example, for the substitution $q_1=\cdots=q_r=q^3$ it is {\underbar{not true}} that the form is preserved. In particular, for any function $\epsilon(q,M_1,\cdots,M_r,L)$ such that $\epsilon(q,M_1,\cdots,M_r,L)\to 1$ as $M_1,\cdots ,M_r,L$ all $\to\infty$ we have\medskip
\begin{align*}
&\sum_{\nu_1, \nu_2, \cdots, \nu_r, k\ge 0}a^{\nu_1}_1a^{\nu_2}_2\cdots a^{\nu_r}_rc^kq^{3(T_{\nu_1}+T_{\nu_2}+\cdots+T_{\nu_r})+T_k+k(\nu_1+\nu_2+\cdots \nu_r)}{{L}\brack {N+k}}_{q}{{N}\brack {\nu_1, \cdots, \nu_r}}_{q^3}\\
&\qquad\qquad\qquad\qquad\times \quad {{N+k}\brack {N,k}}_{q}\epsilon(q,M_1,\cdots,M_r,L)\\
&\not=\sum_{\nu_1, \nu_2, \cdots, \nu_r\ge 0}a^{\nu_1}_1a^{\nu_2}_2\cdots a^{\nu_r}_rq^{3(T_{\nu_1}+T_{\nu_2}+\cdots+T_{\nu_r})}{{M_1}\brack {\nu_1}}_{q^3}{{M_2}\brack {\nu_2}}_{q^3}\cdots {{M_r}\brack {\nu_r}}_{q^3}(-q)_N(-cq^{N+1})_{L-N}. 
\end{align*}\\
This is because, if we take $a_2=\cdots=a_r=c=0$ and let $M_1\to \infty$, the coefficient of $a_1$ in LHS is $\frac{q^3}{1-q}$ and the coefficient of $a_1$ in RHS is $\frac{q^3(1+q)}{1-q^3}$. This is to say, if we set $q_1=\cdots=q_r=q^{m'}$, only the case $m'=2$ gives the form of \eqref{eq71}.\bigskip

\textbf{Remark 8.6:} Like \eqref{eq710}, using the bound shifting method on \eqref{eq817} gives, for any $L\in\mathbb{Z}$, 
\begin{align*}
&\sum_{\nu_1, \cdots, \nu_r, k\geq 0}a^{\nu_1}_1\cdots a^{\nu_r}_rc^kq^{2(\sum_{i=1}^{r}T_{\nu_i})+T_k+kN}{{L}\brack {N+k}}_{q}{{N}\brack {\nu_1, \cdots, \nu_r}}_{q^2} {{N+k}\brack {N,k}}_{q}(-q^{L-N+1})_{N}\\
&=\sum_{\nu_1, \nu_2, \cdots, \nu_r\ge 0}a^{\nu_1}_1a^{\nu_2}_2\cdots a^{\nu_r}_rq^{2(\sum_{i=1}^{r}T_{\nu_i})}\prod_{i=1}^r{{L-\sum_{i< j\leq r}\nu_j}\brack {\nu_i}}_{q^2}(-q)_N(-cq^{N+1})_{L-N}.\numberthis\label{eq818}
\end{align*}

\textbf{Remark 8.7:}
In \cite{YAKA25} we will consider \cite[Eq. (5.5)]{KAGEABG95} as an instance of the more general equation below.
\[
H(\nu_1,\dots,\nu_r,k)=q^{2T_{\nu_1}+\cdots +2T_{\nu_r}+T_k+k(\nu_1+\nu_2+\cdots \nu_r)}{{\nu_1+\cdots+\nu_r}\brack {\nu_1, \cdots, \nu_r}}_{q^2}\times{{\nu_1+\cdots+\nu_r+k}\brack {\nu_1+\cdots+\nu_r,k}}_{q}\]
We will interpret the above function combinatorially in a manner that has similarities and differences to the interpretation in \cite{KAGEABG95}. This interpretation highlights the combintorial significance of the forms \refeq{eq71} and \refeq{eq81b}.\bigskip

\textbf{Work of Berkovich-Uncu:} Berkovich and Uncu \cite{ABAKU19JNT,ABAKU19AC,ABAKU20,ABAKU22}, provided three separate finite polynomial identities that imply Capparelli's identities as limiting cases. The original polynomial identities were found by imposing bounds on the combinatorial constructions that Kur\c{s}ung\"{o}z \cite{KK19} introduced. Their proofs used an automated deduction method and recurrences. The proven polynomial identities also led to the discovery of infinite hierarchies of sum-product identities that include Capparelli's identities \cite{ABAKU19JNT}, \cite{ABAKU20}. Similarly, Uncu \cite{AKU20}, by refining Kur\c{s}ung\"{o}z's construction, proved a polynomial identity that implies Schur's partition theorem. The results and methods in this paper are different from these related works.\bigskip

\textbf{Concluding thought:} Two things have been accomplished in this paper: (i) a unification of Schur's theorem and Lebesgue identity, and (ii) an infinite hierarchy of $q$-hypergeometric identities of which the initial ones are the identities for the partition theorems of Lebesgue and Capparelli. The unification of the Schur and Lebesgue partition theorems was motivated by the fact that the basic idea behind both of their combinatorial proofs was the same - namely, a method of Bressoud suitably adapted to each of the two. Since the combinatorial proof of Capparelli's theorem in \cite{KAGEABG95} also uses the technique of Bressoud, it is our desire to unify the theorems of Schur, Lebesgue, and Capparelli by a common scheme.\bigskip

\textbf{Acknowledgements:} We thank George Andrews, Alexander Berkovich, and Ali Uncu for stimulating discussions and helpful remarks. We also appreciate the help of Aritram Dhar, who converted our amstex file to LaTeX and made pointed remarks on Theorem C5. Additionally, we are grateful to Ole Warnaar for his comments and suggestions. Lastly, we thank both referees for their helpful remarks. In particular, we thank Referee 2 for a careful reading of, and observations on, the statement of Theorem C5. This ultimately put us on the path of refining the statement.\bigskip

\textbf{Conflict of Interest:} We have 
no relevant financial or non-financial interest 
to disclose. Thus, there is no conflict of interest 
or competing interests.

\newpage

\bibliographystyle{amsplain}


\end{document}